\documentstyle[12pt]{article}

\newcommand{\argsup}{\mathop{\rm argsup}\limits}

\newcommand{\const}{\mathop{\rm const}\limits}

\newcommand{\meas}{\mathop{\rm meas}\limits}

\newcommand{\supp}{\mathop{\rm supp}\limits}

\newcommand{\Ent}{\mathop{\rm Ent}\limits}

\begin{document}

\begin{center}

{\bf BOUNDEDNESS OF OPERATORS IN BILATERAL GRAND LEBESGUE SPACES} \\
\vspace{3mm}
{\bf with exact and weakly exact constant calculation.}\\

\vspace{3mm}

 $ {\bf E.Ostrovsky^a, \ \ L.Sirota^b } $ \\

\vspace{4mm}

$ ^a $ Corresponding Author. Department of Mathematics and computer science, Bar-Ilan
University, 84105, Ramat Gan, Israel.\\
\end{center}
E - mail: \ eugostrovsky@list.ru\\
\begin{center}
$ ^b $  Department of Mathematics and computer science. Bar-Ilan University,
84105, Ramat Gan, Israel.\\
\end{center}
E - mail: \ sirota@zahav.net.il\\

\vspace{3mm}

{\bf Abstract}.  In this article we investigate an action of some operators (not
necessary to be linear or sublinear) in the  so-called (Bilateral) Grand Lebesgue
Spaces (GLS), in particular, double weight Fourier operators, maximal operators,
imbedding operators etc. \par
 We intend to calculate an exact or at least weak exact values for correspondent
imbedding constant. \par
We obtain also interpolation theorems for GLS spaces.\par
 We construct several examples to show the exactness of offered estimations.\par
 In two last sections we introduce anisotropic Grand Lebesgue Spaces, obtain some
 estimates for Fourier two-weight inequalities and calculate Boyd's multidimensional
 indices for this spaces.\par

\vspace{3mm}

 \par

\vspace{3mm}

\vspace{3mm}

\section{Introduction. Notations. Problem Statement.} \par

\vspace{3mm}

  Let $ (X, {\cal A},\mu) $ and $ (Y, {\cal B}, \nu) $  be a two measurable spaces
with sigma-finite non - trivial measures $ \mu, \nu. $ For the measurable real valued
functions $ f(x), \ x \in X, f: X \to R; \ g(y), \ y \in Y, g: Y \to R $
the symbols $ |f|_p = |f|_p(X,\mu) $ and correspondingly $ |g|_q = |g|_q(Y,\nu) $
will denote the usually $ L_p $  and $ L_q $ norms:
$$
|f|_p = |f|L_p(X, \mu) = \left[ \int_X |f(x)|^p \ \mu(dx) \right]^{1/p}, \ p \ge 1.
\eqno(1.1a)
$$

$$
|g|_q = |g|L_q(Y, \nu) = \left[ \int_Y |g(y)|^q \ \nu(dy) \right]^{1/q}, \ q \ge 1.
\eqno(1.1b)
$$
 Let $ U $ be an operator, not necessary to be linear or sublinear, defined on arbitrary
measurable simple bounded function $ f: X \to R $ with values in the set of measurable
functions $ g: Y \to R. $ We impose the following important condition on the operator
$ U. $ \par
{\bf Condition A, "moment inequality".} There exist an intervals
$ (a,b), 1 \le a < b \le \infty,  \ (c,d), 1 \le c < d \le \infty, $ and strictly
monotonically continuous function $  q = q(p),  q: (a,b) \to (c,d), $  for which the
operator $ U $ is bounded as the operator from the space $ L_p(X) = L_p(X,\mu) $ onto the
space $ L_q(Y) = L_q(Y,\nu):$

$$
|Uf|_{q(p)} \le K(p) \ |f|_p, \ \forall p \in (a,b) \ \Rightarrow K(p) < \infty. \eqno(A)
$$

{\it We assume hereafter that the interval $ (a,b) $ is the {\it maximal interval} for
 which the inequality (A) there holds. \par
 We denote the inverse function to the function $ p \to q(p) $ as $ r = r(q); \
 r: (c,d) \to (a,b), q(r(p)) = p. $ \par
 Further, we understand as the constant $ K(p) $  its {\it minimal value}, namely: }
$  K(p) = |U| \left(L_p \to L_{q(p)} \right) $ {\it or equally}

$$
K(p) = \sup_{f \in L_p(X), f \ne 0} \left[\frac{|Uf|_{q(p)}}{ |f|_p} \right] =
|U|_{p \to q(p)}.
$$

 There are many examples of operators satisfying the condition (inequality)
 (A) with calculated the exact (or exactly evaluated)
  value of constant $ K(p): $ classical Hardy-Littlewood
 inequalities \cite{Hardy1}, with its  Bradley's generalizations \cite{Bradley1};
 integral \cite{Okikiolu1}, \cite{Dunford1}, chapter 5, section 11, p. 567-580;
 in particular, convolution  operators (regular and
 singular) \cite{Beckner1}, oscillating  integrals \cite{Stein3}, potentials 
 \cite{Jorgens1}, p. 270-271;  Sobolev's imbedding
 operators with modern generalizations \cite{Maz'ya1}, \cite{Mitrinovich1}; 
 Poincare-Sobolev's inequalities  \cite{Badiale1}, \cite{Ghoussoub1}, \cite{Vassilev1},
 maximal operators with Muckenhoupt's \cite{Muckenhoupt1} generalization,  Hilbert's 
 transforms \cite{Bennet1}, chapter 3, Fourier's transform operators (discrete and 
 continuous) \cite{Boas1}, \cite{Chen1}, \cite{Chen2}, \cite{Titchmarsh1},
\cite{Yu1}; pseudodifferential (Fourier integral) operators \cite{Taylor1} etc. \par

 For instance, let us consider the following integral operator:

 $$
 U_o[f](x) = \int_R |t|^{\mu-1} h(x \cdot t) f(t) dt,
 $$

$$
1 \le p \le q,\ 1/q = \mu - 1/p, 1/\sigma := 1+\mu - 2/p.
$$
 Here

 $$
 K(p) = \left[ \int_R |t|^{-\sigma(p-1)/p} |h(t)|^{\sigma} \ dt \right]^{1/\sigma},
 $$
see \cite{Okikiolu1}, p. 222.\par
 Let us concern  the so-called potential operators, of a Riesz's type (non-homogeneous):
 $$
 P[f](x) = \int_{R^n} a(x) \ b(y) \ |x-y|^{\beta-n} \ f(y) \ dy.
 $$
 Let 
 
 $$
 \frac{p-1}{p} + \frac{1}{p_0} < 1, \ \frac{1}{q} + \frac{1}{q_0} < 1, \ \beta \in (0,n),
 $$
 
 $$
 \frac{1}{p} = \frac{1}{q} + \frac{1}{p_0} \frac{1}{q_0} - \frac{\beta}{n}.
 $$
 It is known that 
 
 $$
 |P|(L_p \to L_q) \le C \cdot |a|_{p_0} \cdot |b|_{q_0},
 $$
 see \cite{Jorgens1}.\par
  For the following modification of Fourier transform:

 $$
 \hat{f}(y) = \int_{R^n} e^{-2\pi i x y}f(x) dx
 $$
W.Beckner in \cite{Beckner1} proved the following estimates with exact constants
computations:

$$
|\hat{f}|_{p_1} \le A(p) \ |f|_p, \ p \in (1,2], \ p_1= p/(p-1),
$$
where

$$
A(p) \stackrel{def}{=} \left[ \frac{p^{1/p}}{p_1^{1/p_1}}  \right]^{n/2};
$$

$$
|f*g|_r \le \left( A(p) A(q) A(r_1)  \right)^n \ |f|_p \ |g|_q, 1/r = 1/p +1/q - 1.
$$
S.K.Pichorides in \cite{Pichorides1} proved the following inequality for the Hilbert's
transform $ H[f] $ with sharp value of constant:

$$
|H[f]|_p \le \Lambda(p) \ |f|_p, \ 1 < p < \infty,
$$

$$
\Lambda(p) = \tan(\pi/(2p)), 1 < p \le 2; \ \Lambda(p) = \cot(\pi/(2p)), 2 < p < \infty.
$$

 Analogous estimates  are known for Sobolev's (Hardy-Littlewood-Sobolev) imbedding theorems,
see \cite{Talenti1}.  We refer here only the so-called {\it fractional} Sobolev's inequality
\cite{Hichem1}, \cite{Hongwei1}:

$$
|f|_q \le K_S(p) \cdot \left|  \left[\sqrt{-\Delta}\right]^s [f] \   \right|_p,
$$

$$
\left[\sqrt{-\Delta}\right]^s [f](y) \stackrel{def}{=}
F^{-1} \left(|x|^s \ F[f](x) \right)(y),
$$

$$
K_S(p) = \pi^{s/2} \cdot \frac{\Gamma((n-s)/2)}{\Gamma((n+s)/2)} \cdot
\left\{ \frac{\Gamma(s)}{\Gamma(n/2)} \right\}^{s/n},
$$

$$
0<s<n, \ 1 < p < n/s, \ q = \frac{pn}{n-sp}, \ \Delta[f](x_1,x_2, \ldots,x_n)=
\sum_{m=1}^n \frac{\partial^2}{\partial x_m^2}f.
$$

\vspace{3mm}

 Next example: pseudodifferential operators, on the other terms: Fourier integral
 operators. We refer to the book of Taylor M.E. \cite{Taylor2}, article \cite{Nagase1}.\par
  Consider for example  the pseudodifferential operator $ P(D) $ of a view

 $$
 P(D) [f](x) = \int_{R^n} e^{ix\xi} p(\xi) F[f](\xi) \ d \xi.
 $$
 If

 $$
 \forall k: |k| \le \Ent(n/2) + 1 \ \Rightarrow
 \sup_{R>0} \int_{R < \xi < 2R} | \ |\xi|^{|k|} D^k P(\xi) \ |^2 d \xi < \infty,
 $$
 $ \Ent(z) $ denotes the integer part of a variable $ z, $ then
 $ P(D): L_p \to L_p, \ p \in (1,\infty) $  and

 $$
 |P(D)|(L_p \to L_p) \le C \cdot \frac{p^2}{p-1}.
 $$

Too modern results are obtained in  \cite{Ferreira1}.\par

\vspace{3mm}

Other example: weight inequalities for Fourier transform; we follow B.Muckenhoupt
\cite{Muckenhoupt1}. \par
 Denote as usually

 $$
 F[f](y) = (2\pi)^{-n/2} \int_{R^n} e^{i x y} \ f(x) \ dx.
 $$
Let $ u = u(y), \ v = v(x) $ be two positive integrable functions (weights).  We cite
here the inequality of a view:

$$
\left[\int_{R^n} |F[f](y)|^q \ u(y) \ dy   \right]^{1/q} \le K_M(p,q) \
\left[\int_{R^n} |f(x)|^p \ v(x) \ dx  \right]^{1/p}, \ p,q  \in (1,\infty).
$$
 {\it Suppose  here that  $ p \le q; $ the alternative case will be studied
 further.} \par
  For instance, let $ p \in (1,2]; $ we denote for some positive constants $ A, B $

  $$
  I(A,B)= \sup_{r > 0} \left\{ \left[\int_{u(y) > B \ r} u(y) \ dy \right] \cdot
\left[ \int_{v(x) < A \ r^{p-1}} v(x)^{-1/(p-1)} \ dx  \right]\right\}.
  $$

 Muckenhoupt in \cite{Muckenhoupt1} proved in particular that for $ q = p/(p-1) $

$$
\left[\int_{R^n} |F[f](y)|^q \ u(y) \ dy   \right]^{1/q} \le  \frac{CI(A,B)}{p-1}
\cdot \left[\int_{R^n} |f(x)|^p \ v(x) \ dx  \right]^{1/p}, \ p,q  \in (1,\infty).
$$
 This inequality may be rewritten as follows.

$$
\left[\int_{R^n} |F[f](y)|^q \ u(y) \ dy   \right]^{1/q} \le \inf_{A,B > 0}
\frac{CI(A,B)}{p-1}
\cdot \left[\int_{R^n} |f(x)|^p \ v(x) \ dx  \right]^{1/p}, \ p,q  \in (1,\infty).
$$
 For the classical Calderon-Zygmund singular integral
 operators $ U $ (and its commutators with vector
 fields) are known the asymptotically exact up to multiplicative constants of a view

$$
|U|(L_p \to L_p) \le  C \ \frac{p^2}{p-1}, \ p \in (1,\infty).
$$
see \cite{Otriz-Caraballo1}, \cite{Lerner1}. \par

 Analogous estimates for different modifications of Fourier integral operators
for example for the {\bf maximal} Fourier operator

 $$
 |\sup_N S_N[f]|_p \le C \frac{p^4}{(p-1)^2} \ |f|_p, \ p \in (1,\infty),
 $$
where $ S_N[f](x) $ denotes the $ N \ - $ th partial Fourier sum for the function
$ f: [0, 2 \pi] \to R $  see, e.g. in the book  of Reyna  \cite{Reyna1}. \par

\vspace{3mm}
 The last considered in this section
 example belongs to E.M.Stein and G.Weiss \cite{Stein2}. Consider the following
 singular integral operator $ S $ with homogeneous kernel $  K = K(x,y), \ x,y \in R^n: $

 $$
 S[f](x) = \int_{R^n} K(x,y) f(y) dy,
 $$
where $ K(x,y) \ge 0, $

$$
\forall \lambda > 0 \ \Rightarrow K(\lambda x, \lambda y) = \lambda^{-n}
K(x,y),
$$

$$
\forall Y \in SO(n) \ \Rightarrow K(Yx,Yy) = K(x,y).
$$

 Proposition:

$$
|S|(L_p \to L_p) = \int_{R^n} K(x,e_1) |x|^{-n/p^/} dx,  \ e_1 = (1,0,0,\ldots,0).
$$
 This result generalized classical inequalities of Hardy-Littlewood. \par

\vspace{3mm}

 {\bf  Our  aim is extrapolation of the "moment" inequality (A) on the so-called
 Grand Lebesgue Spaces (GLS) with exact or at least weak exact constants  calculations.} \par

\vspace{3mm}

 The papier is organized as follows. In the next section we recall used facts
about Grand Lebesgue Spaces, formulate and prove the main result of this paper. \par
 In the third section we obtain the GLS boundedness of the so-called maximal transform.
 The fourth section is devoted to the weight Fourier operators boundedness in GLS spaces. \par
 In the next section we consider the interpolation inequalities for operators acting
 in GLS spaces. We based here on the classical interpolation results belonging to
 Riesz-Thorin and Marcinkiewicz.\par
  In the sixth section we generalize obtained results on the "multidimensional" case and
 prove the exactness of the relation between parameters in the
Lebesgue spaces Fourier and convolution weight operators inequalities, by means of a
so-called "dilation method". The next section is devoted to the
generalization of obtained results on the anisotropic spaces; we prove also the exactness
of our estimates.\par
 The object of the eight section is generalization of  previous results on the case when
the weight functions are arbitrary  continuous regular varying.\par
The last section contains some concluding remarks.\par

\vspace{3mm}

\section{Grand Lebesgue Spaces. Main result.}

\vspace{3mm}

  Let $ (X,A,\mu) $ be again measurable space with sigma-finite
non - trivial measure $ \mu. $ We recall in this section  for readers conventions some
definitions and facts from the theory of GLS spaces.\par
\vspace{2mm}

Recently, see \cite{Fiorenza1}, \cite{Fiorenza2}, \cite{Fiorenza3}, \cite{Iwaniec1},
\cite{Iwaniec2}, \cite{Kozachenko1},\cite{Liflyand1}, \cite{Ostrovsky1}, \cite{Ostrovsky2}
 etc.  appear  so-called Grand Lebesgue Spaces $ GLS = G(\psi) =G\psi =
G(\psi; A,B), \ A,B = \const, A \ge 1, A < B \le \infty, $ spaces consisting
on all the measurable functions $ f: X \to R $ with finite norms

$$
     ||f||G(\psi) \stackrel{def}{=} \sup_{p \in (A,B)} \left[ |f|_p /\psi(p) \right].
     \eqno(2.1)
$$
Here $ \psi(\cdot) $ is some continuous positive on the {\it open} interval
$ (A,B) $ function such that

$$
  \inf_{p \in (A,B)} \psi(p) > 0, \ \psi(p) = \infty, \ p \notin (A,B). \eqno(2.2)
$$
We will denote
$$
 \supp (\psi) \stackrel{def}{=} (A,B) = \{p: \psi(p) < \infty, \} \eqno(2.3)
$$

The set of all $ \psi $  functions with support $ \supp (\psi)= (A,B) $ will be
denoted by $ \Psi(A,B). $ \par
  This spaces are rearrangement invariant, see \cite{Bennet1}, and
are used, for example, in the theory of probability  \cite{Kozachenko1},
\cite{Ostrovsky1}, \cite{Ostrovsky2}; theory of Partial Differential Equations
\cite{Fiorenza2}, \cite{Iwaniec2};  functional analysis \cite{Fiorenza3},
\cite{Iwaniec1},  \cite{Liflyand1}, \cite{Ostrovsky2}, \cite{Ostrovsky16},
\cite{Talagrand1}; theory of Fourier series \cite{Ostrovsky1}, theory of martingales
\cite{Ostrovsky2},  mathematical statistics \cite{Sirota1}, \cite{Sirota2}, \cite{Sirota3},
\cite{Sirota4}, \cite{Sirota5}, \cite{Sirota6}, \cite{Sirota7}, \cite{Sirota8}; theory of
approximation \cite{Ostrovsky7} etc.\par
 Notice that in the case when $ \psi(\cdot) \in \Psi(A,B),  $ a function
 $ p \to p \cdot \log \psi(p) $ is convex, and  $ B = \infty, $ then the space
$ G\psi $ coincides with some {\it exponential} Orlicz space. \par
 Conversely, if $ B < \infty, $ then the space $ G\psi(A,B) $ does  not coincides with
 the classical rearrangement invariant spaces: Orlicz, Lorentz, Marzinkievicz etc.\par

\vspace{2mm}

{\bf Remark 2.0.} If we define the {\it degenerate } $ \psi_r(p), r = \const \ge 1 $
function as follows:
$$
\psi_r(p) = \infty, \ p \ne r; \psi_r(r) = 1
$$
and agree $ C/\infty = 0, C = \const > 0, $ then the $ G\psi_r(\cdot) $ space coincides
with the classical Lebesgue space $ L_r. $ \par

\vspace{2mm}

 Let $ \xi: X \to R $ be measurable function such that for some constants
 $ (A,B), 1 \le A < B \le \infty  \ |\xi|_p < \infty. $  The {\it natural } function
 for the function $ \xi = \xi(x) \  \psi_{\xi}(p) $  may be defined by formula

 $$
 \psi_{\xi}(p) := |\xi|_p, \ p \in (A,B).
 $$

 Analogously, let $ \xi = \xi(t,x), \ t \in T, \ T $  is arbitrary set,  be a
 {\it family} of a measurable functions such that
 $$
 \exists (A,B): 1 \le A < B \le \infty, \ \forall p \in (A,B) \Rightarrow
   \sup_{t \in T} \ |\xi(t,\cdot)|_p < \infty.
 $$
The {\it natural } function $ \xi = \xi(x) \  \psi_{\xi}(p) $ for the {\it family }
$ \xi(\cdot) $ denotes by definition

$$
\psi_{\xi}(p) := \sup_{t \in T} |\xi(t,\cdot)|_p, \ p \in (A,B).
$$

\vspace{3mm}

{\sc Hereafter we will denote by $ c_k = c_k(\cdot), C_k = C_k(\cdot), k = 1,2,\ldots, $
with or without subscript, some positive finite non-essentially "constructive"
constants, non necessarily at the same at each appearance .\par
 We will denote also by the symbols $ K_j = K_j(n,p,q,\ldots) $  essentially
positive finite functions depending only on the} $ n,p,q,\ldots. $ \par

\vspace{2mm}

 Let  $ \psi = \psi(p), \ p \in (a,b) $ be some function from the set $ G\Psi(a,b); $
we define a new function

$$
\psi_1(q) = K(r(q)) \times  \psi(r(q)).  \eqno(2.4)
$$

{\bf Theorem 2.1.} {\it Let} $ f \in G\psi; $ {\it then}

$$
||U[f]||G\psi_1 \le 1 \times ||f||G\psi, \eqno(2.5)
$$
{\it and the constant "1"  in the inequality (2.5) is the best possible. }\par
{\bf Proof of the upper estimate}  is very simple. Let $ f \in G\psi;  $  we can
suppose $  ||f||G\psi = 1. $ It follows from the direct definition of the
norm in the GLS that

$$
\forall p \in (a,b) \ \Rightarrow \ |f|_{p,X} \le \psi(p).
$$
 We obtain from the condition (A)

 $$
 |U[f]|{q(p)} \le K(p) \psi(p),
 $$
or equally

$$
|U[f]|_q \le K(r(q)) \psi(r(q)) = \psi_1(q) = \psi_1(q) ||f||G\psi,
$$

$$
||U[f]||G\psi_1 = \sup_{q \in (c,d)} [|U[f]|_q/\psi_1(q)] \le ||f||G\psi.
$$

{\bf Proof of the lower estimate}. We denote

$$
V(\psi,f) = \left[ \frac{||U[f]||G\psi_1}{||f||G\psi} \right],
\overline{V} = \sup_{\psi \in G\psi(a,b)} \sup_{f \in G\psi} V(\psi,f),
$$
and define as usually

$$
||U||(G\psi \to G\psi_1) = ||U|| = \sup_{f \in G\psi, f \ne 0}
\left[\frac{||U[f]||G\psi_1}{||f||G\psi} \right],
$$
then

$$
||U|| = \sup_{f \in G\psi, f \ne 0}
\left[\frac{\sup_p|U[f]|_{q(p)}/\psi_1(q(p))}{\sup_p|f|_p/\psi(p)} \right].
$$

 The proposition or theorem (2.1) may be rewritten as follows: $ \overline{V} = 1; $
we know that $ \overline{V} \le 1; $ it remains to prove an opposite inequality. \par
 Let us choose

$$
q_0 = \argsup_{q \in (c,d)} \sup_{f \in G\psi, f \ne 0}
\left[ \frac{|U[f]|_q}{K(r(q)) |f|_{r(q)}} \right],
$$

 $$
 f_0 = \argsup_{f \in G\psi} \left[\frac{|U[f]|_{q_0}}{K(r(q_0)) |f|_{r(q_0)}} \right],
 \ \psi_0(p) = |f_0|_p = \psi_{f_0}(p),
$$

$$
\psi_{0,1}(q) = K(r(q)) \times  \psi_0(r(q)).
$$

 It follows from the definition of the function $ K = K(p) $ that

 $$
 |U[f_0]|_{q_0} = K(r(q_0)) |f_0|_{r(q_0)},
 $$
and we can suppose without loss of generality  that the function $ f_0 $ here exists.\par
 Further,

$$
\overline{V} = \sup_{\psi \in G\psi(a,b)} \sup_{f \in G\psi} V(\psi,f) =
\sup_{\psi \in G\psi(a,b)} \sup_{f \in G\psi}
\left[ \frac{||U[f]||G\psi_1}{||f||G\psi} \right] =
$$

$$
\sup_{\psi \in G\psi(a,b)} \sup_{f \in G\psi}
\left[ \frac{\sup_q|U[f]|_q /\psi_1(q)}{ \sup_p |f|_p /\psi(p)} \right] \ge
 \frac{\sup_q|U[f_0]|_q /\psi_{0,1}(q)}{ \sup_p |f_0|_p /\psi_0(p)} \ge
$$
$$
\frac{|U[f_0]|_{q_0}}{\psi_{0,1}(q_0)}  =
\frac{K(r(q_0)) |f_0|_{q_0}}{K(r(q_0))|f_0|_{q_0} } = 1, \eqno(2.6)
$$
as long as $ f_0 \ne 0. $ \par

 In the case when the function $  f_0 $ does'nt exist, we conclude that for arbitrary
 $ \epsilon \in (0,1/2) $ there exists a {\it family } of a measurable functions
 $ f_{\epsilon} = f_{\epsilon}(x) $  for which

 $$
 |U[f_{\epsilon}]|_{q_{\epsilon}} > (1-\epsilon) \ K(r(q_{\epsilon})) \
 |f_{\epsilon}|_{r(q_{\epsilon})},
 $$
and therefore

$$
\overline{V} \ge 1-\epsilon.
$$

{\bf Remark 2.1.} In all known examples the value $ p_0 $ tends to the boundaries
of the set $ (a,b); $ may be $ p_0 \to \infty $ if $ b = \infty. $ \par
{\bf Remark 2.2.} There are many cases when $ q = p $ and following
$ (c,d) = (a,b), $ for example Hilbert transform or maximal Fourier operators etc.
 Obviously, here

 $$
 \psi_1(p) = K(p) \ \psi(p),  \ p \in (a,b).
 $$

{\bf Remark 2.3.} If we know instead the exact value of constant $ K(p) $ only the
{\it weak inequality } of a view

$$
0 \le c_1 \le \sup_{p \in (a,b)}  \frac{|U|_{p \to q(p)}}{K(p)} \le C_2 < \infty, \eqno(2.7)
$$
then evidently

$$
c_1 \le ||U||(G\psi \to G\psi_1) \le C_2, \eqno(2.8)
$$
(the "weak value" of  constant.) \par
{\bf Remark 2.4.} It may be investigated analogously the case of inequality

$$
|U_1[f]|_{q(p)} \le K(p) \ |U_2[f]|_p, \ \forall p \in (a,b) \ \Rightarrow K(p) <
\infty, \eqno(B)
$$
where $ U_1, U_2 $ are some operators and the second operator $ U_2 $ is not invertible.
 For example, the inequality (B) is  true for Sobolev's imbedding theorems. \par
{\bf Remark 2.5.} We consider here the case when the value $ q = q(p) $ is'nt unique.
More exactly, let for some interval $ (a,b), 1 \le a < b \le \infty $ of the values
$ p $ there exists a non-trivial interval $ Q(p) = (Q_1(p), Q_2(p)), 1 \le Q_1(p)
 < Q_2(p) \le \infty $ of the values $ q $  and positive finite function
 $ K_Q(p,q) $ such that

$$
|U[f]|_q \le K_Q(p,q) \ |f|_p. \eqno(2.9)
$$
 Let for some non-zero function $ \psi \in \Psi(a,b) \ f \in G\psi, $ then

 $$
 |U[f]|_q \le K_Q(p,q) \ \psi(p) \ ||f||G\psi, \eqno(2.10)
 $$
 and let $ \nu = \nu(q) $  be another function from the set $ \Psi(c,d). $  We obtain
from (2.10) after dividing on $ \nu(q) \cdot ||f||G\psi: $

$$
\frac{|U[f]|_q}{\nu(q) \ ||f||G\psi} \le \frac{K(p,q) \ \psi(p)}{\nu(q)}.
$$
 We get taking supremum over $ q: $
$$
\frac{||U[f]||G\nu}{||f||G\psi} \le \sup_q \left[\frac{K(p,q) \ \psi(p)}{\nu(q)} \right].
\eqno(2.11)
$$
 As long as the inequality is true for each values $ p, $ we conclude eventually:

 $$
 ||U[f]||G\nu \le
 \inf_{p \in (a,b)} \sup_{q \in Q(p)} \left[\frac{K(p,q) \ \psi(p)}{\nu(q)} \right]
 \cdot ||f||G\psi =:\overline{K} \ ||f||G\psi, \eqno(2.12)
 $$
if obviously $ \overline{K} < \infty. $\par
{\bf Remark 2.6.} Let us consider a particular case of inequality (2.9); namely, suppose
for any constants $ A,B > 0 $
$$
|U[f]|_q \le C \cdot A^{1/q} \ B^{-1/p} \  |f|_p, \ p \in (a,b), q \in (c,d).\eqno(2.13)
$$
 We get repeating at the same considerations:

 $$
\frac{|U[f]|_q \cdot B^{1/p}}{\nu(q)\psi(p)} \le C \frac{A^{1/q}|f|_p}{\nu(q)\psi(p)}.
\eqno(2.14)
 $$
 Recall that the {\it fundamental function} of the space $ G\psi $ has a view:

 $$
 \phi(G\psi, \delta) = \sup_{p \in (a,b)} \left[ \frac{\delta^{1/p}}{\psi(p)} \right].
 $$
We have taking supremum of the inequality (2.14) over $ q $ and $ p: $

$$
\frac{||U[f]||G\nu}{\phi(G\nu,A)} \le C \cdot \frac{||f||G\psi}{\phi(G\psi,B)}. \eqno(2.15)
$$

 The inequality (2.15) was used in the approximation theory, see \cite{Ostrovsky7}.\par
{\bf Remark 2.7.} The lower bound in the theorem 2.1 may be obtained very simple if we
allow to consider the degenerate $  \psi \ - $ functions. \par
 {\bf Example of non-uniqueness. } \par
 We return to the article of Muckenhoupt \cite{Muckenhoupt1}, but we represent more
 general case:  $ q \ne p/(p-1). $ Namely, let here $ 1 < p < 2 $  and
 $ p \le q < p_1 = p/(p-1). $ Denote by $ u^*(y) $ the non-increasing rearrangement
 modification of the function $  u(y), $  and define the following integral:

 $$
 J(A,B) = \sup_{r > 0} \left[J_1(B,r) J_2(A,r) \right],
 $$
 where
 $$
 J_1(B,r)= \int_{ \{ [y u^*(y) ]^{p_1/q} > B r x  \} } u^*(y) \ dy,
 $$

$$
J_2(A,r) =  \left[ \int_{ \{v(x)< A r^{p-1} \} } v(x)^{-1/(p-1)} \ dx \right]^{q/p_1}.
$$

 Suppose $ J(A,B) < \infty $ for some constants $ A,B > 0 $ (Muckenhoupt's condition).
 We conclude after some computations from the article  \cite{Muckenhoupt1}:

 $$
\left[ \int_{R^n} u(y) \ |F[f](y)|^q \ dy  \right]^{1/q} \le K_{MM}(p,q) \cdot
\left[\int_{R^n} v(x) \ |f(x)|^p \ dx  \right]^{1/p},  \eqno(2.16)
 $$
where

$$
K_{MM}(p,q) = \frac{C(A,B;J(A,B))}{(p-1)(p_1-q) }.\eqno(2.17)
$$

\vspace{2mm}

\section{ Boundedness of maximal operator in GLS}

\vspace{2mm}

 Let $ U = U_{\theta} = U_{\theta}[f], \ \theta \in \Theta, $ where
 $ \Theta = \{\theta\} \ $  is arbitrary set, be a family of linear
 (or at least  sublinear) operators.  The sublinear (in general case)  operator of a view

 $$
 \overline{U}[f](x) = \sup_{\theta \in \Theta} |U_{\theta}[f](x)|, \ x \in X \eqno(3.1)
 $$
will be called {\it maximal} operator for the family $ U_{\theta}, $ if is correctly
defined on some Banach functional space of a (measurable) functions $ f: X \to R. $ \par
 We consider in this section the boundedness of operator $ \overline{T}[f] $ in
 Grand Lebesgue Spaces.\par
 Note that the case of maximal Fourier transform and some other singular operators is
 observed, e.g., in \cite{Fiorenza3}, \cite{Ostrovsky21}. \par
 We  consider only the case of classical maximal {\it centered ball} Hardy -
 Littlewood  operator in the ordinary Euclidean space $ X =R^d: $

 $$
{\cal M}[f](x) = \sup_{B:|B| \in (0,\infty)} \left[ |B|^{-1}\int_B |f(x)| dx \right],
\ \eqno(3.2)
 $$
 where "supremum"  is calculated over all Euclidean balls $ B $ with center at the
 point $ x;  \ |B| := \meas(B) \in (0,\infty).  $\par
 Let $ \psi \in G\psi(1,\infty); $ we define

 $$
 \psi^{(1)}(p) \stackrel{def}{=} \psi(p) \cdot \frac{p}{p-1}, \ p \in (1,\infty). \eqno(3.3)
 $$
 Note that $ \psi^{(1)}(\cdot) \in G\psi(1,\infty). $ \par
 {\bf Theorem 3.1.}  If $ f \in G\psi, $ then

 $$
 || {\cal M}[f] ||G\psi^{(1)} \le C_3 ||f||G\psi,
 $$
for some {\it absolute}, i.e.  not depending on the dimension $ d $ weakly exact
non-trivial constant $ C_3. $ \par
{\bf Proof.} We will use theorem 2.1, more exactly remark 2.3. The classical result
belonging to  E.M.Stein E.M and J.O.Str\"omberg \cite{Stein1} states that the constant
$ C_2 $ in (2.7) for considered maximal operator $ {\cal M} $ is bounded over all the
dimensions $ d: \ \sup_d C_2(d) =:C_3 < \infty. $ \par
 In detail, E.M.Stein E.M and J.O.Str\"omberg \cite{Stein1} proved:

$$
 | {\cal M} |_{p \to p} \asymp \frac{p}{p-1}, \ p \in (1,\infty].
$$

The lower estimate for the constant $ c_1 $ in (2.7) is also uniform over $ d $ bounded
from below, which may be proved by consideration of an example

$$
f_0 = I(||x|| \le 1), \ ||x|| = \sqrt{\sum_{i=1}^d x^2_i},
$$
$$
I(x \in A) = I_A(x) =1, \ x \in A; \ I(x \in A) = I_A(x) = 0, x \notin A.
$$

We have for the function $ f_0: $

$$
| {\cal M}[f_0] |_p \ge c_1 \frac{p}{p-1} |f_0|_p.
$$
 This completes the proof of theorem 3.1.\par

\vspace{3mm}
\section{ Weight Pitt - Beckner - Okikiolu  inequalities for GLS spaces}
\vspace{3mm}

Let $ x = \vec{x}  \in R^n $ be $ n- \ $ dimensional vector,
$ n = 1,2,\ldots. $ which consists on the $ l, l \ge 1 $ subvectors
$ \{ x_j, \} \ j = 1,2,\ldots,l: $

$$
x = (x_1,x_2, \ldots, x_l), \  x_j = \vec{x_j} \in R^{m_j},\  \dim(x_j) = m_j \ge 1.
\eqno(4.1)
$$
 Let also $ \alpha = \vec{\alpha} = \{ \alpha_1,\alpha_2, \ldots,\alpha_l \} $ and
$ \beta = \vec{\beta} = \{ \beta_1,\beta_2, \ldots,\beta_l \}$  be two fixed $ l \ - $
dimensional vectors such that

$$
\alpha_j \in [0,1), \ \beta_j \in [0,1),\ \alpha_j + \beta_j \le 1, \ j = 1,2,\ldots,l.
$$
 We denote as ordinary by $ |x_j| $ the Euclidean norm of the vector $ x_j;$

 $$
 |x|^{-\alpha} = \prod_{j=i}^l |x_j|^{-\alpha_j}
 $$
 and analogously for the vector $ y \in R^n $

 $$
 |y|^{\beta} = \prod_{j=i}^l |y_j|^{\beta_j}. \eqno(4.2)
 $$

 Obviously,
 $$
 \sum_{j=1}^l m_j = l, \ |x| = (\sum_{j=1}^l |x_j|^2)^{1/2}.
 $$
 We define alike to the book Okikiolu \cite{Okikiolu1}, p. 313-314, see also
 \cite{Beckner2}  the so-called {\it double weight Fourier} transform
 $ F_{\alpha,\beta}[f](x) $ by the following way:

$$
F_{\alpha,\beta}[f](x)= (2\pi)^{-n/2} |x|^{\alpha} \int_{R^n} |y|^{\beta} \
f(y) \ e^{i x y} \ dy, \eqno(4.3)
$$

$$
F[f](x)= (2 \pi)^{-n/2} \int_{R^n} f(y) \ e^{i x y} \ dy = F_{0,0}[f](x),
$$
where $ xy $ denotes the inner product of the vectors $ x,y: xy = \sum_{k=1}^n x_k y_k. $
\par

 The inequality of a view (relative Lebesgue measures in whole space $ R^n $ )

 $$
 |F_{\alpha,\beta}[f]|_q \le K_{PBO}(p) |f|_p \eqno(4.4a),
 $$
 or equally

 $$
 | \ |y|^{-\alpha} \ F[f](y)|_q \le K_{PBO}(p) \ | \ |x|^{\beta} \ f(x)|_p, \eqno(4.4)
 $$
is said to be (generalized) weight Pitt - Beckner - Okikiolu (PBO) inequality. \par
 We will understood as customary as the value  $ K_{PBO}(p)$ its minimal value:

 $$
 K_{PBO}(p) =
 \sup_{f \ne 0, f \in L_p(R^n)} \left[ \frac{|F_{\alpha,\beta}[f]|_q }{|f|_p} \right].
  \eqno(4.5).
 $$

  As before, $ p \in (p_0,p_1), \ p_0,p_1 = \const, 1 \le p_1 < p_2 \le \infty, \
 q \in (q_0,q_1), \ 1 \le q_0 < q_2 \le \infty $ and $ q = q(p) $ is some uniquely defined
 continuous strictly monotonic function (if there exists). \par
 There are many publications about this inequality, see also, for instance,
\cite{Beckner2}, \cite{Boas1}, \cite{Chen1}, \cite{Chen2}, \cite{Gorbachov1},
\cite{Gord1}, \cite{Jurkat1}, \cite{Leindler1}, \cite{Liflyand1}, \cite{Liflyand2},
\cite{Sagher1}, \cite{Yu1} etc.\par
 A very interest application of the PBO inequality in the quantum mechanic are described
in the articles \cite{Beckner1}, \cite{Herbst1}.\par
{\it We  consider in this section only "one-dimensional" case, i.e. when} $ l=1. $ \par
The general case $ l \ge 2 $ will be considered further. \par
 Let us introduce the following important conditions in the domain
 $$
 p_0 = n/(n-\beta), p_1 = \infty,  \ q_0 = 1,  q_1 = n/\alpha:
 $$
 $$
 \beta-\alpha = n - n \left(\frac{1}{p} + \frac{1}{q} \right), \eqno(4.6)
 $$

 $$
 p > p_0, \ 1 \le q < q_0,
 $$

which defined uniquely the continuous function $ q = q(p). $ \par
 The condition (4.6) was imposed by many authors, see \cite{Beckner1}, \cite{Gorbachov1},
\cite{Liflyand2} etc. \par
\vspace{3mm}
{\bf Theorem 4.1.}  \par
\vspace{3mm}
{\bf A.} The conditions 4.6 is necessary for PBO inequality (4.4).\par
{\bf B.} If  the condition 4.6 is satisfied, then for $ p>p_0  $

$$
C_1(\alpha,\beta,n)
\left[\frac{p}{p-p_0} \right]^{(\alpha+\beta)/n} \le K_{PBO}(p)\le
$$

$$
C_2(\alpha,\beta,n) \left[\frac{p}{p-p_0} \right]^{\max(1,(\alpha+\beta)/n)}. \eqno(4.7)
$$

{\bf C.} Both the boundaries in the inequality (4.7), lower and upper, are attainable. \par
\vspace{3mm}

{\bf Proof. \ Part B.} \par
{\bf 0.} The upper bound for the coefficient $ K_{PBO}(p) $ follows from
the direct computation in the article of W.Beckner \cite{Beckner2}. For instance, in the
case when $ \alpha = \beta $ or correspondingly $ q = p/(p-1) $ W.Beckner computed the
exact value of this constant.\par
 {\bf 1.} In order to obtain the lower bound in (4.7), we consider the following example.
 Let us introduce the following subsets (generalized truncated "segments") of whole space
$ R^n = \{ x = \vec{x}: \} $

 $$
 D(c_0,c_1,c_2) = \{x: |x_j| \ge c_0, \ |x_j|/|x| \in [c_1,c_2]  \}, \ j = 1,2,\ldots,n;
 \eqno(4.8)
 $$

 $$
 G(c_3,c_4,c_5) = \{y: 0 < |y_j| \le c_3, \ |y_j|/|y| \in [c_4,c_5]  \}, \ j = 1,2,\ldots,n;
 \eqno(4.9)
 $$

$$
0 < c_0 < \infty; \ 0 < c_1 < 1 < c_2 < \infty; \
0 < c_3 < \infty; \ 0 < c_4 < 1 < c_5 < \infty;
$$

$$
D = D(c_1,c_2) = D(1,c_1,c_2); \  G = G(c_4,c_5) = G(1,c_4,c_5).
$$

{\bf 2.} We consider the following example.  Let $ f_0 = f_0(x), \ x \in R^n $ be
even function such that

$$
f_0(x) = f_{0,n}(x) = \frac{I_D(x)}{\prod_{j=1}^n |x_j| }. \eqno(4.10)
$$
We have using multidimensional polar (spherical) coordinates:

$$
| \ |x|^{\beta} \ f_0(x)|_p^p \le \int_{D(c_1,c_2)} |x|^{p(\beta-n)} dx  \le C_3(n)
\int_{C_4(n)}^{\infty} r^{n-1+p(\beta-n)} dr, \eqno(4.11)
$$
as long as inside the domain $ D(c_1,c_2) $ is true the following inequality:

$$
\prod_{j=1}^n |x_j| \ge C_5(n,c_1,c_2) \ |x|^n.
$$
 We conclude taking the integral (4.11) that if $ p > n/(n-\beta)  $ then

 $$
| \ |x|^{\beta} \ f_0(x)|_p \le C_6(n,\beta) (p-n/(n-\beta))^{-1/p} \le
C_7(n,\beta) \left( p/(p-p_0)  \right)^{-(n-\beta)/n}.\eqno(4.12)
 $$

{\bf 3.} We investigate in this pilcrow  the behavior of the
Fourier transform of the function $ f_0(y) $ as $ |y| \to 0, \ |y| \le 1: $

$$
F[f_0](y) = \int_D e^{ixy} \frac{dx}{\prod_{j=1}^n |x_j|} =
C_8(n) \int_{D} \cos(xy) \ [\prod_{j=1}^n I(x_j > 0)] \ \frac{dx}{\prod_{j=1}^n |x_j|}.
$$
 We restrict ourself only the case when $ y \in G(c_3,c_4,c_5).$ \par
 We obtain after the substitution $ x_j = v_j/y_j: $

$$
F[f_0](y)  = C_9 \int_{G(y)}  \cos(\sum_j v_j) \ \frac{dv}{\prod_{j=1}^n |v_j|} \sim
C_9 \int_{G(y)} \frac{dv}{\prod_{j=1}^n |v_j|},
$$
where $ G(y) $ is an image of the set  $ D $ under our substitution.\par
 Since for all the values $ y, \ y \in G(c_3,c_4,c_5) $ the domain $  G(y) $ contains
 the set $ D(c_6|y|,c_7, c_8), $ we have for the values $ y \in G(c_3,c_4,c_5) $ using
 again the multidimensional polar (spherical) coordinates as before:

$$
F[f_0](y) \ge C_{10} \int_{D(c_6|y|,c_7, c_8)} \frac{dv}{\prod_{j=1}^n |v_j|} \ge
$$

$$
C_{11} \int_{C_{12}|y|}^{C_{13}} \rho^{n-1 - n} \ d \rho \ge C_{14} |\log |y| \ |.
\eqno(4.13)
$$

{\bf 4.} We estimate here the left-hand side of inequality (2.7) for our operator. Namely,

$$
| \ |y|^{-\alpha} \ F[f](y) |^q_q \ge \int_{0 < |y| \le c_7} C_{15}^q |y|^{-\alpha q} \
|\log|y| \ |^q dy \ge
$$

$$
C_{16}^q \int_{0}^{C_{17}} z^{n-1-\alpha q} \ |\log z|^q dz \ge
C_{18}^q(n,\alpha) (n-\alpha q)^{-q-1},
$$
or equally

$$
| \ |y|^{-\alpha} \ F[f_0](y) |_q \ge C_{19}(n,\alpha) \ (q_0-q)^{-1-1/q} \ge
$$

$$
C_{20}(n,\alpha,\beta) (p-p_0)^{-1-1/q} \ge C_{21}
\left[\frac{p}{p-p_0}\right]^{-1-\alpha/n}; \eqno(4.14)
$$
here $ q \in [1, n/\alpha) $ and correspondingly $ p \in (n/(n-\beta),\infty). $ \par
{\bf 5.} We conclude after dividing $ | \ |y|^{-\alpha} \ F[f_0](y) |_q $
over  $ | \ |x|^{\beta} \ f(x)|_p:  $

$$
K_{PBO}(p) \ge \frac{ | \ |y|^{-\alpha} \ F[f_0](y) |_q}{ | \ |x|^{\beta} \ f_0(x)|_p} \ge
C_{21}(n,\alpha,\beta) \left[\frac{p}{p-p_0} \right]^{-\beta/n -\alpha/n}.\eqno(4.15)
$$
 This completes the proof of the proposition {\bf B } of theorem 4.1.\par

\vspace{3mm}

{\bf Proof of the assertion  C.} \par
{\bf Lower bound.}  The lower bound in the inequality (4.7) is attained as
$ \alpha \to 0+, \ \beta \to 0+, $ for instance, in the Hausdorff-Young inequality

$$
|F[f]|_{p/(p-1)} \le C^n |f|_p, \ p \in (1,2],
$$
where $ C $ is an absolute constant (calculated by W. Beckner \cite{Beckner1},
\cite{Beckner2}); here $ \alpha = \beta =  0. $\par
 The {\bf upper} bound is attained, e.g., in the case when $ \alpha = \beta > 0,
 q = p/(p-1), $ see \cite{Beckner2}.\par
  At the same upper estimate may be obtained from the classical Wiener-Paley theorem:

$$
\int_R |y|^{p-2} |F[f](y)|^p \ dy \le K_{WP}^p(p) \ \int_R |f(x)|^p \ dx,
$$
see  \cite{Krein1}, chapter 5, section 5; in this case

$$
K_{WP}(p) \asymp \frac{p}{p-1}, \ p \in (1,2].
$$

 \vspace{3mm}

{\bf Remark 4.1.} Notice that our counter-example is'nt a radial function! \par
\vspace{3mm}

 As a {\bf corollary:} let $ \psi(\cdot) \in G\psi(p_0,\infty); $ we define a new function

$$
\psi_{PBO}(p) = \frac{p \psi(p)}{p-p_0}.
$$
 We assert:

 $$
 ||F_{\alpha,\beta}[f]||G\psi_{PBO} \le C_2(\alpha,\beta,n) ||f||G\psi, \eqno(4.16)
 $$
and the last inequality is weak logarithmicaly  exact.\par

\vspace{3mm}

{\bf Proof of  the part A.} \par
 We will use (here and further) the so-called "dilation method",  introduced, e.g. by
G.Talenti \cite{Talenti1} for a finding of an optimal constant in the Sobolev's
imbedding theorems. Indeed, let $ f, f:R^n \to R $ be some non-zero function from the
Schwartz space $ S(R^n) $, for which the inequality PBO (4.4) there holds. We define
the {\it family} $ T_{\lambda}[\cdot], \ \lambda \in (0,\infty) $  of linear
dilations operators of a view

$$
T_{\lambda}[f](x) = f(\lambda x). \eqno(4.17)
$$
 Note that if $ f \in  S(R^n), $ then $ \forall \lambda \in (0,\infty) \
 \Rightarrow T_{\lambda}[f] \in S(R^n). $ \par
 We have:
 $$
 F[T_{\lambda} \ f] = \lambda^{-n} F[f], \eqno(4.18)
 $$

$$
| \ |x|^{\beta} \ T_{\lambda} \ f |_p =  \lambda^{-n/p - \beta} | \ |x|^{\beta} \ f |_p,
\eqno(4.19),
$$

$$
| \ |y|^{-\alpha} \ F[T_{\lambda} \ f]|_q = \lambda^{-n + n/q - \alpha}
| \ |y|^{-\alpha}\ f]|_q . \eqno(4.20)
$$
 We get after substituting onto (4.4) for the function $ T_{\lambda} \ f: $

$$
\lambda^{-n + n/q - \alpha} | \ |y|^{-\alpha} \ F[f]|_q  \le K_{PBO}(p)
 \lambda^{-n/p - \beta} | \ |x|^{\beta} \ f |_p, \eqno(4.21).
$$
 As long as the last inequality may be satisfied for all the values
$ \lambda \in (0,\infty) $ only in the case when

$$
-n + n/q - \alpha = -n/p - \beta,
$$
or equally
$$
\beta - \alpha = n -n(1/p + 1/q),
$$
Q.E.D.

\vspace{3mm}

\section{ Interpolation of operators in Grand Lebesgue spaces}

\vspace{2mm}

  Let $ (X, {\cal A},\mu) $ and $ (Y, {\cal B}, \nu) $  be again two measurable spaces
with sigma-finite non-trivial measures $ \mu, \nu; $ we assume in addition that both
the measures $ \mu, \nu $ are resonant in the sense described, e.g. in monograph Bennet and
Sharpley \cite{Bennet1}, chapter 1,2.  Recall that  the notion "resonant" of the
 measure $ \mu $ means either  atomlessness this measure or discreteness $ \mu $ such that
 all the atoms have equal measures. \par

\vspace{3mm}

{\bf A. Classical version.}\par

\vspace{3mm}

  We consider first of all the interpolation of linear operators in the spirit of the
 classical theorem belonging to Riesz  and Thorin.
   Let $ p_0,p_1,q_0,q_1 $ be fixed numbers such that $ 1 \le p_0,p_1,q_0, q_1 < \infty,
   \ p_0 \le q_0, p_1 \le q_1. $ \par
We can and will suppose without loss of generality that $ q_0 < q_1, \ p_0< p_1. $ \par
 We consider in this section a linear operator from the set of all measurable functions
 $ f: X \to R $ into the set of measurable functions  $ f: Y \to R $ such that

   $$
   |T[f]|_{p_0} \le M_0 |f|_{q_0}, \ |T[f]|_{p_1} \le M_1 |f|_{q_1}, M_0,M_1 = \const.
   \eqno(5.1)
   $$
 On the other words, the operator $ T $ is simultaneously of strong-type $ (p_0,q_0) $
and $ (p_1,q_1). $ \par
 Define a functions for  $ q \in (q_0,q_1): $

$$
r_{RT}(q) = \frac{p_0 p_1 (q_1-q_0) q}{p_0 q_1(q-q_0) + p_1 q_0 (q_1-q) },
$$

$$
\Theta(q) = \frac{q_1(q-q_0)}{q(q_1-q_0)},
$$
so that

$$
1-\Theta(q) = \frac{q_0(q_1-q)}{q(q_1-q_0)};
$$

$$
M_{\theta}(q) = 2 M_0^{1-\Theta(q)} \ M_1^{\Theta(q)}.
$$

{\bf Theorem 5.1.} Let $ \psi \in G\psi(p_0,p_1); $  we define a new function

$$
\psi_{RT} = \psi_{RT}(q) = M_{\theta}(q)\cdot \psi \left(r_{RT}(q) \right).\eqno(5.2)
$$
Proposition: for each function from the space $ G\psi(p_0,p_1) $ there holds:

$$
||T[f]||G\psi_{RT} \le ||f||G\psi.\eqno(5.3)
$$

{\bf Proof} follows from theorem 2.1 and from the classical theorem of Riesz and Thorin
\cite{Bennet1}, p. 185:

$$
|T[f]|_q \le M_{\theta}(q) \cdot |f|_p, \ p \in (p_0,p_1), \ q \in (q_0,q_1).
$$
where

$$
\frac{1}{p} = \frac{1-\theta}{p_0} + \frac{\theta}{p_1}, \
\frac{1}{q} = \frac{1-\theta}{q_0} + \frac{\theta}{q_1}. \eqno(5.4)
$$
 {\bf Remark 5.1.} This result is weakly exact, e.g., when $ T $ is Fourier transform
on the whole  real line (or whole space $ R^n.) $ Theorem of Hardy - Young asserts that

 $$
 |T[f]|_q \le C \ |f|_p, \ p \in (1,2),   \ q = p/(p-1).
 $$

\vspace{3mm}

{\bf Lorentz version.}

\vspace{3mm}

We recall first of all the definition of the Lorentz spaces over the triple (say)
$ (X, \cal{A}, \mu). $ Let the numbers $ p,q $ be a given such that $ 1 \le p \le \infty,
0 < q \le \infty. $ \par
 We denote as customary for each measurable function $ f:X \to R $ by $ f^*(t),
 \ t \in (0,\infty), $ the left inverse function  to the distribution  function
 for its absolute value:
 $$
 f^*(t) = h^{-1}(t), \ h(t) = \mu \{x \in X, \ |f(x)| > t \}.
 $$
 The space $ L_{p,q} = L_{p,q}(X) $ consists, by definition, see \cite{Bennet1}, p.
213-214, on all the functions $ f, f:X \to R $  with finite norm

$$
|f|_{p,q} = \left[  \int_0^{\infty} \left[t^{1/p} f^*(t)  \right]^q dt/t \right]^{1/q},
\eqno(5.5)
$$
if $ q  < \infty, $ and

$$
|f|_{p,\infty} = \sup_{t > 0}  t^{1/p} f^*(t).
$$

 Recall that the linear operator $ T $ is said to be of weak type $ (p,q), $ if

 $$
 |T[f]|_{q,\infty} \le M |f|_{p,1}, \ M = \const < \infty.
 $$

{\bf Theorem 5.2.} Let $ p_0,p_1,q_0,q_1 $  be fixed number such that
$  1 \le p_0 < p_1 < \infty, \ 1 \le q_0 < q_1 \le \infty.  $ Let also $ T $ be
simultaneously  weak type $ (p_0,q_0) $ and $ (p_1,q_1): $
$$
|T[f]|_{q_i,\infty} \le M_i |f|_{p_i,1}, \ i=0,1. \eqno(5.6)
$$

We introduce for each function $ \psi \in G\psi(p_0,p_1)$ a new function

$$
\psi_M(q) = \psi \left(r_{RT}(q) \right) \ [(q-q_0)(q_1-q)]^{-1}, \
q \in (q_0,q_1).
$$
 Proposition:

$$
||T[f]||G\psi_M \le C(p_0,p_1,q_0,q_1) \ \max(M_0,M_1) \ ||f||G\psi.\eqno(5.7)
$$

{\bf Proof} is at the same as the proof of theorem 5.1. We will use only instead the
classical theorem of Riesz - Thorin the interpolation theorem belonging to
Marcinkiewicz  (see \cite{Bennet1}, p. 226). Namely, if the conditions (5.6) are
satisfied, then

$$
|T[f]|_q \le \frac{c(p_0,p_1,q_0,q_1)}{\theta(1-\theta)} \ \max(M_0,M_1) \ |f|_p,
$$
where as before $ p \in (p_0,p_1), \ q \in (q_0,q_1), $

$$
\frac{1}{p} = \frac{1-\theta}{p_0} + \frac{\theta}{p_1}, \
\frac{1}{q} = \frac{1-\theta}{q_0} + \frac{\theta}{q_1}.
$$

{\bf Remark 5.2}. The assertion of theorem (5.2) is in general case non-improvable, see,
e.g. \cite{Ostrovsky21},  \cite{Ostrovsky27}.\par

\vspace{3mm}

\section{Multidimensional case. "Dilation method".}

\vspace{3mm}
 We consider in this section the multidimensional (vector): $ l \ge 2 $
 generalization of PBO inequality (4.4) in the sense of section 4. Namely, we
 investigate the inequality of a view

 $$
 | \ |y|^{- \vec{\alpha}} F[f](y)|_q \le K_{n,\vec{\alpha},\vec{\beta}}(p) \
 | \ |x|^{\vec{\beta}} \ f(x)|_p,
 $$
 or for simplicity

 $$
 | \ |y|^{- \alpha} F[f](y)|_q \le K_{n,\alpha,\beta}(p) \
 | \ |x|^{\beta} \ f(x)|_p. \eqno(6.1)
 $$

 Recall that $ x = \vec{x}  \in R^n $ be $ n- \ $ dimensional vector,
$ n = 1,2,\ldots; $ which consists on the $ l, l \ge 1 $ subvectors
$ \{ x_j, \} \ j = 1,2,\ldots,l: $

$$
x = (x_1,x_2, \ldots, x_l), \  x_j = \vec{x_j} \in R^{m_j},\  \dim(x_j) = m_j \ge 1;
$$
$$
 \alpha = \vec{\alpha} = \{ \alpha_1,\alpha_2, \ldots,\alpha_l \}
$$
and
$ \beta = \vec{\beta} = \{ \beta_1,\beta_2, \ldots,\beta_l \}$  be two fixed $ l \ - $
dimensional vectors such that

$$
\alpha_j \in [0,1), \ \beta_j \in [0,1), \ j = 1,2,\ldots,l, \ \alpha_j + \beta_j \le 1.
$$
 We denote as ordinary by $ |x_j| $ the Euclidean norm of the vector $ x_j;$

 $$
 |x|^{-\alpha}:= \prod_{j=i}^l |x_j|^{-\alpha_j}
 $$
 and analogously for the vector $ y \in R^n $

 $$
 |y|^{\beta}:= \prod_{j=i}^l |y_j|^{\beta_j}.
 $$

\vspace{3mm}
{\bf Theorem 6.1.} \par
 {\bf 1.} The inequality (6.1) may be satisfied for each function $ f \in S(R^n) $ only
when

$$
\frac{\beta_j - \alpha_j}{m_j} = \const, \ j \in [1,2,\ldots,l], \eqno(6.2)
$$

 {\bf 2.} If the condition (6.2) is satisfied, then

 $$
  \beta_j - \alpha_j = m_j \left(1 - \frac{1}{p} - \frac{1}{q} \right), \eqno(6.3)
 $$

 $$
 C_1(\alpha,\beta,n)  \prod_{j=1}^l
 \left[ \frac{p}{(p-m_j/(m_j-\beta_j))}\right]^{(\alpha_j+\beta_j)/m_j}  \le
  K_{n,\alpha,\beta}(p) \le
 $$

$$
C_2(\alpha,\beta,n) \prod_{j=1}^l
\left[\frac{p}{(p-m_j/(m_j-\beta_j))} \right]^{\max(1,(\alpha_j+\beta_j)/m_j)}, \eqno(6.4)
$$

 $$
 p > \max_{j=1,2,\ldots,l} \frac{m_j}{m_j-\beta_j}.
 $$

{\bf Proof } of the second proposition is at the same as in the one-dimensional case
(theorem 4.1); the upper bound for the coefficient $ K_{n,\alpha,\beta}(p) $ may be
obtained analogously considerations of W.Beckner \cite{Beckner1}. \par
 The example $ \hat{f}(x) $ for the lower bound may be constructed as a product of the
 $ m_j $ dimensional examples:

$$
\hat{f}(x) = \prod_{j=1}^l f_{0,m_j}(x_j).
$$
 It remains to prove the necessity of the condition (6.3).

 We will use a {\it multidimensional} generalization of dilation method. \par
 It is enough to consider only the two-dimensional case: $ l=2. $
 The inequality (6.1) then has a view:

 $$
 | \ |y_1|^{- \alpha_1} \ |y_2|^{-\alpha_2} \ F[f](y)|_q \le K_{n,\alpha,\beta}(p) \
 | \ |x_1|^{\beta_1} \ |x_2|^{\beta_2} f(x)|_p. \eqno(6.5)
 $$

 Let $ \lambda_1, \lambda_2 $ be two arbitrary independent positive numbers.
We define the  multidimensional dilation operator $ T_{\lambda_1, \lambda_2}[f] $
 as follows:

 $$
 T_{\lambda_1, \lambda_2}[f](x_1,x_2) = f(\lambda_1 x_1, \lambda_2 x_2).\eqno(6.6a)
 $$
 We have:

 $$
 F[T_{\lambda_1, \lambda_2}[f]](y_1,y_2) =
 \lambda_1^{-m_1} \lambda_2^{-m_2} \
 F[f](y_1/\lambda_1,y_2/\lambda_2);\eqno(6.6b)
 $$

$$
| \ |y_1|^{-\alpha_1} \ |y_2|^{-\alpha_2} \ F[T_{\lambda_1, \lambda_2}[f]](y_1,y_2)|_q =
$$

$$
\lambda_1^{m_1/q -m_1 - \alpha_1} \ \lambda_2^{m_2/q - m_2 - \alpha_2}
| \ |y_1|^{-\alpha_1} \ |y_2|^{-\alpha_2} \ F[f(y_1,y_2)]|_q . \eqno(6.6c)
$$

$$
| \ |x_1|^{\beta_1} \ |x_2|^{\beta_2} T_{\lambda_1, \lambda_2}[f](x_1,x_2)|_p =
$$
$$
\lambda_1^{-\beta_1 - m_1/p} \ \lambda_2^{-\beta_2 - m_2/p}
| \ |x_1|^{\beta_1} \ |x_2|^{\beta_2} \ f(x_1,x_2)|_p; \eqno(6.6d).
$$

 We obtain substituting into (6.5) the function $ T_{\lambda_1, \lambda_2}[f] $ instead
the function $ f: $

$$
\lambda_1^{m_1/q -m_1 - \alpha_1} \ \lambda_2^{m_2/q - m_2 - \alpha_2}
| \ |y_1|^{-\alpha_1} \ |y_2|^{-\alpha_2} \ F[f(y_1,y_2)]|_q \le
$$

$$
\lambda_1^{-\beta_1 - m_1/p} \ \lambda_2^{-\beta_2 - m_2/p}
| \ |x_1|^{\beta_1} \ |x_2|^{\beta_2} \ f(x_1,x_2)|_p. \eqno(6.7)
$$

Since the values $ \lambda_1, \ \lambda_2 $ are arbitrary positive, the inequality
may be satisfied only when

$$
m_1/q - m_1 - \alpha_1 = -\beta_1 - m_1/p,
$$
$$
 m_2/q - m_2 - \alpha_2 = -\beta_2 - m_2/p,
$$
or equally

$$
(\beta_1 - \alpha_1)/m_1 = 1-1/p-1/q = (\beta_2 - \alpha_2)/m_2.\eqno(6.8)
$$
This completes the proof of theorem 6.1.\par
{\bf Remark 6.1.} At the same result as theorem 6.1, up to simple chaining of
variables and functions, may be obtained from the
considerations Okikiolu \cite{Okikiolu1}, p. 313-318. The  condition (6.2) is there
presumed (in another notations). \par

\vspace{3mm}

\section{PBO inequalities for anisotropic Lebesgue spaces. }

\vspace{3mm}

We consider in this section the case of some generalization of PBO inequality
when the condition (6.2) is'nt satisfied. \par

 In order to formulate and prove the generalization of the inequality (6.1),
we will use the so-called {\it anisotropic} Lebesgue spaces.  More detail information
about this spaces see in the books \cite{Besov1}, chapter 16,17;  \cite{Leoni1}, chapter
 11; using for us theory of operators interpolation in this spaces see in \cite{Besov1},
chapter 17,18.\par

 Let $ \vec{p} = (p_1,p_2,\ldots,p_l) $ be $ l \ - $ dimensional vector such that

 $$
 p_j > m_j/(m_j-\beta_j) =:p_{0,j}, \ q_j \in [1, m_j/\alpha_j) =:q_{0,j}, \
 j=1,2,\ldots,l.
 $$
 We denote the set all the values $ \vec{p} $ as
 $ Q = Q \left(p_{0,1}, p_{0,2}, \ldots, p_{0,l} \right). $ \par

  Let also $ u = u(x), \ x \in R^n $ be measurable function: $ u:R^n \to R.$ Recall
 that the anisotropic Lebesgue space $ L_{\vec{p}} $ consists on all the functions
 $ f $ with finite norm

 $$
 |f|_{\vec{p}}\stackrel{def}{=}
\left( \int_{R^m_1} dx_1 \left( \int_{R^{m_2}}dx_2 \ldots  \left( \int_{R^{m_l}}
|f(\vec{x})|^{p_1} dx_1 \right)^{p_2/p_1} \right)^{p_3/p_2} \ldots \right)^{1/p_l}.
\eqno(7.1)
 $$
 Note that in general case

 $$
 |f|_{p_1,p_2} \ne |f|_{p_2,p_1},
 $$
but
$$
|f|_{p,p} = |f|_p.
$$
 Observe also that if $ f(x_1,x_2) = g_1(x_1)\cdot g_2(x_2) $ (condition of factorization),
 then

 $$
 |f|_{p_1,p_2} = |g_1|_{p_1} \cdot |g_2|_{p_2},
 $$
(formula of factorization). \par

 Let $ \psi = \psi(\vec{p}) $ be some continuous positive on the set
 $ Q = Q \left(p_{0,1}, p_{0,2}, \ldots, p_{0,l} \right) $ function such that

$$
  \inf_{p \in Q} \psi(p) > 0, \ \psi(p) = \infty, \ p \notin Q. \eqno(7.2)
$$
We denote the set all of such a functions as $ \Psi_Q. $ \par
 The (multidimensional, anisotropic)  Grand Lebesgue Spaces $ GLS = G_Q(\psi) =G_Q\psi $
space consists on all the measurable functions $ f: R^n \to R $ with finite norms

$$
||f||G_Q(\psi) \stackrel{def}{=} \sup_{\vec{p} \in Q}
\left[ |f|_{\vec{p}} /\psi(\vec{p}) \right].
     \eqno(7.3)
$$

\vspace{3mm}

 The object of our investigation in this section is an inequality of a view

$$
 | \ |y|^{- \vec{\alpha}} F[f](y)|_{\vec{q}} \le
 K_{n,\vec{\alpha},\vec{\beta}}(\vec{p}) \
 | \ |x|^{\vec{\beta}} \ f(x)|_{\vec{p}},
 $$
 briefly:
  $$
 | \ |y|^{- \alpha} F[f](y)|_{\vec{q}} \le K_{n,\alpha,\beta}(\vec{p}) \
 | \ |x|^{\beta} \ f(x)|_{\vec{p}}, \eqno(7.4)
 $$
(the classical version), or correspondingly

 $$
 || \ |y|^{- \alpha} F[f](y) \ ||G\nu_Q \le  \
 || \ |x|^{\beta} \ f(x)  \ ||G\psi_Q, \eqno(7.5)
 $$
(Grand Lebesgue Spaces version). \par
 We introduce before the formulation the following notations. The connection between
 $ p_j  $ and $ q_j, \ = 1,2,\ldots, l $ may be described by the formulae

 $$
 \beta_j - \alpha_j = m_j \left(1-\frac{1}{p_j} - \frac{1}{q_j} \right).\eqno(7.6)
 $$

\vspace{3mm}
{\bf Theorem 7.1.}\par
{\bf A.} The condition (7.6) is necessary and sufficient for the existing and finiteness
of the constant $ K_{n,\alpha,\beta}(\vec{p}) $ for the inequality (7.5).\par
{\bf B.} If the condition (7.6) is satisfied, and

$$
\alpha_j, \beta_j \ge 0, \ \alpha_j + \beta_j \le 1,
$$
then the sharp (minimal) value of the
coefficient $ K_{n,\alpha,\beta}(\vec{p}) $ satisfies the inequalities

 $$
 C_1(\vec{\alpha},\vec{\beta,\vec{m}})  \prod_{j=1}^l
 \left[ \frac{p_j}{(p_j-m_j/(m_j-\beta_j))}\right]^{(\alpha_j+\beta_j)/m_j}  \le
  K_{n,\alpha,\beta}(\vec{p}) \le
 $$

$$
C_2(\vec{\alpha},\vec{\beta},\vec{m}) \prod_{j=1}^l
\left[\frac{p_j}{(p_j-m_j/(m_j-\beta_j))} \right]^{\max(1,(\alpha_j+\beta_j)/m_j)},
\eqno(7.7)
$$

$$
\ \vec{p} \in Q = Q \left(p_{0,1}, p_{0,2}, \ldots, p_{0,l} \right).
$$

\vspace{3mm}
{\bf Proof.} We will use again the {\it multidimensional} generalization of dilation
method. \par
 It is enough to consider only the two-dimensional case: $ l=2. $
 The inequality (7.4) has then a view:

 $$
 | \ |y_1|^{- \alpha_1} \ |y_2|^{-\alpha_2} \ F[f](y)|_{q_1,q_2} \le K_{n,\alpha,\beta}(p) \
 | \ |x_1|^{\beta_1} \ |x_2|^{\beta_2} f(x)|_{p_1,p_2}. \eqno(7.8)
 $$

 Let $ \lambda_1, \lambda_2 $ be two arbitrary independent positive numbers.
We define the  multidimensional dilation operator $ T_{\lambda_1, \lambda_2}[f] $
 as follows:

 $$
 T_{\lambda_1, \lambda_2}[f](x_1,x_2) = f(\lambda_1 x_1, \lambda_2 x_2).\eqno(7.9a)
 $$
 We have as before:

 $$
 F[T_{\lambda_1, \lambda_2}[f]](y_1,y_2) =
 \lambda_1^{-m_1} \lambda_2^{-m_2} \
 F[f](y_1/\lambda_1,y_2/\lambda_2);\eqno(7.9b)
 $$

$$
| \ |y_1|^{-\alpha_1} \ |y_2|^{-\alpha_2} \
F[T_{\lambda_1, \lambda_2}[f]](y_1,y_2)|_{q_1,q_2} =
$$

$$
\lambda_1^{m_1/q_1 -m_1 - \alpha_1} \ \lambda_2^{m_2/q_2 - m_2 - \alpha_2}
| \ |y_1|^{-\alpha_1} \ |y_2|^{-\alpha_2} \ F[f](y_1,y_2)|_{q_1,q_2}. \eqno(7.9c)
$$

Further,

$$
| \ |x_1|^{\beta_1} \ |x_2|^{\beta_2} T_{\lambda_1, \lambda_2}[f](x_1,x_2)|_{p_1,p_2} =
$$
$$
\lambda_1^{-\beta_1 - m_1/p_1} \ \lambda_2^{-\beta_2 - m_2/p_2}
| \ |x_1|^{\beta_1} \ |x_2|^{\beta_2} \ f(x_1,x_2)|_{p_1,p_2}; \eqno(7.9d).
$$

 We obtain substituting into (7.4) the function $ T_{\lambda_1, \lambda_2}[f] $ instead
the function $ f: $

$$
\lambda_1^{m_1/q_1 -m_1 - \alpha_1} \ \lambda_2^{m_2/q_2 - m_2 - \alpha_2}
| \ |y_1|^{-\alpha_1} \ |y_2|^{-\alpha_2} \ F[f](y_1,y_2)|_{q_1,q_2} \le
$$

$$
\lambda_1^{-\beta_1 - m_1/p_1} \ \lambda_2^{-\beta_2 - m_2/p_2}
| \ |x_1|^{\beta_1} \ |x_2|^{\beta_2} \ f(x_1,x_2)|_{p_1,p_2}. \eqno(7.10)
$$

Since the values $ \lambda_1, \ \lambda_2 $ are arbitrary positive, the inequality
may be satisfied only when

$$
m_1/q_1 - m_1 - \alpha_1 = -\beta_1 - m_1/p_1,
$$
$$
 m_2/q_2 - m_2 - \alpha_2 = -\beta_2 - m_2/p_2,
$$
or equally

$$
\beta_1 - \alpha_1 = m_1(1-1/p_1-1/q_1),
$$
$$
\beta_2 - \alpha_2 = m_2(1-1/p_2-1/q_2). \eqno(7.11)
$$

 The example $ \hat{f}(x) $ for the lower bound may be constructed as before by means of
factorization property of the anisotropic  norm  as a product of the $ m_j $ dimensional
examples:

$$
\hat{f}(x) = \prod_{j=1}^l f_{0,m_j}(x_j).
$$

This completes the proof of theorem 7.1.\par
\vspace{2mm}
 As a consequence: \par
 {\bf Proposition 7.1.} Let $ \psi_Q \in \Psi_Q; $ we define a new function
 $ \nu_Q = \nu_Q(\vec{p}) \in \Psi_Q, \ \vec{p} \in Q $ as follows:

 $$
 \nu_Q(\vec{p}) = \psi_Q(\vec{p}) \cdot \prod_{j=1}^l
 \left[\frac{p_j}{(p_j-m_j/(m_j-\beta_j))} \right]^{\max(1,(\alpha_j+\beta_j)/m_j)}.
  \eqno(7.12)
$$
 We assert:
$$
 || \ |y|^{- \alpha} F[f](y) \ ||G\nu_Q \le  C(\vec{\alpha},\vec{\beta},\vec{m}) \cdot
 || \ |x|^{\beta} \ f(x)  \ ||G\psi_Q, \eqno(7.13)
 $$

\vspace{3mm}
\section{PBO inequality with regular varying weight.}

\vspace{3mm}

Let $  L =L(z), \ M = M(z), \ z \in (0,\infty)$ be two
{\it slowly varying}
simultaneously as $ z \to 0  $  and as $ z \to \infty $ continuous positive functions:

$$
\lim_{\lambda \to 0} \frac{L(\lambda z)}{L(\lambda)} =
\lim_{\lambda \to \infty} \frac{L(\lambda z)}{L(\lambda)} = 1, \eqno(8.1a)
$$

$$
\lim_{\lambda \to 0} \frac{M(\lambda z)}{M(\lambda)} =
\lim_{\lambda \to \infty} \frac{M(\lambda z)}{M(\lambda)} = 1. \eqno(8.1b)
$$
 We refer reader to the book of Seneta  \cite{Seneta1} to the using further facts
 about regular and slowly varying functions.\par
We investigate in this section the slight generalization of PBO inequality of a view:

 $$
 | \ |y|^{-\alpha} \ M(|y|) \ F[f](y)|_q \le K_{LM}(p) \
 | \ |x|^{\beta} \ L(|x|) \ f(x)|_p,  \eqno(8.2)
 $$
so that both the functions
$$
 z \to |z|^{\beta} \ L(|z|), \ z \to  |z|^{-\alpha} \ M(|z|)
$$
are regular varying simultaneously as $ z \to 0  $  and as $ z \to \infty. $ \par
 For example,

 $$
 L(z) = \max \left( |\log z|^{\theta_1},1 \right), \ M(z) = \max( |\log z|^{\theta_2},1),
 \theta_{1,2} = \const \ge 0.
 $$
\vspace{3mm}
{\bf Theorem 8.1.} \par
\vspace{3mm}
{\bf A.}  The inequality (8.2) is true if and only if

$$
\alpha_j,\beta_j \ge 0, \ \alpha < n, \  \beta < n,
$$

$$
 q < q_0 = n/\alpha, \ p > p_0 = n/(n-\beta),  \eqno(8.3a)
$$

$$
\beta - \alpha = n(1-1/p - 1/q),  \ L(1/\lambda) \asymp  M(\lambda), \lambda \in (0,\infty),
\eqno(8.3b)
$$
in the sense that

$$
0 < \inf_{\lambda > 0} \frac{L(1/\lambda)}{M(\lambda)} \le
\sup _{\lambda > 0} \frac{L(1/\lambda)}{M(\lambda)} < \infty. \eqno(8.3c)
$$

{\bf B.} If the conditions (8.3a), (8.3b), (8.3c) are satisfied, then the function
$ K_{LM}(p) $ satisfies at the same restrictions as the function $ K_{PBO}(p): $ \par

$$
C_3(L,M;\alpha,\beta,n)
\left[\frac{p}{p-p_0} \right]^{(\alpha+\beta)/n} \le K_{LM}(p)\le
$$

$$
C_4(L,M;\alpha,\beta,n) \left[\frac{p}{p-p_0} \right]^{\max(1,(\alpha+\beta)/n) }.
\eqno(8.4)
$$

\vspace{2mm}
{\bf Proof.} The estimates  (8.4) are proved as before, as in \cite{Beckner2}, with as
the same counterexample. The using in \cite{Beckner2} weight interpolation inequalities
see, e.g. in the book \cite{Bennet1}, chapter 4. Namely, if

$$
|U[f]\cdot v_i|_{q_i} \le M_i |f \cdot u_i|_{p_i}, i-0,1;
$$

$$
\frac{1}{p} =  \frac{1-\theta}{p_0} +  \frac{\theta}{p_1},
$$

$$
\frac{1}{q} =  \frac{1-\theta}{q_0} +  \frac{\theta}{q_1},
$$

$$
u=u_0^{1-\theta}\cdot u_1^{\theta}, \  v=v_0^{1-\theta}\cdot v_1^{\theta}, \
u_i,v_i \ge 0,
$$
then

$$
|U[f]\cdot v|_q \le M_0^{1-\theta} \ M_1^{\theta} \ |f \cdot u|_p.
$$
 Note only that if $ v_0, \ v_1 $ are slowly varying, then  the function $ v(\cdot) $
 is also slowly varying.\par
 It remains to prove the assertion {\bf A.} We will use as before the dilation method.\par

 Suppose for sone function $ f \in S(R^n), \ f \ne 0 $

$$
| \ |y|^{-\alpha} \ F[f](y) \ L(|y|) |_q \le K_{LM}(p) \
| \ |x|^{\beta} \ f(x) \ M(|x|)|_p. \eqno(8.5)
$$
 We obtain applying (8.5) to the function $ T_{\lambda} f(x) = f(\lambda x):$

 $$
 \lambda^{-n + n/q - \alpha } | \ |y|^{-\alpha} \ F[f](y) L(|y|/\lambda) \ |_q \le
 K_{LM}(p) \
 \lambda^{-n/p-\beta} | \ |x|^{\beta} \ f(x) M(x \lambda) \ |_p. \eqno(8.6)
 $$
 As long as the functions $ L,M $ are slowly varying, we conclude that as
 $ \lambda \to \infty $ or $ \lambda \to 0+: $

$$
 \lambda^{-n + n/q - \alpha }\ L(1/\lambda) | \ |y|^{-\alpha} \  F[f](y) \ |_q \le
$$
$$
 K_{LM}(p) \ \lambda^{-n/p - \beta} \ M(\lambda) \
  | \ |x|^{\beta} \ f(x) \ |_p. \eqno(8.7)
$$

 The last inequality (8.7) may be satisfied only in the case when

 $$
 -n + n/q - \alpha = -n/p - \beta, \  L(1/\lambda) \asymp M(\lambda),
 $$
Q.E.D.\par
 Analogously may be proved the multidimensional generalization of inequality (8.2).
Indeed, let us consider the following inequality:

 $$
 | \ |y|^{- \vec{\alpha}} F[f](y) \ \prod_{j=1}^l M_j(|y_j|)|_{\vec{q}}
 \le K_{L,M;n,\vec{\alpha},\vec{\beta}}(\vec{p}) \ \cdot
 $$

 $$
 | \ |x|^{\vec{\beta}} \ f(x) \ \prod_{j=1}^l L_j(|x_j|)|_{ \vec{p}}, \eqno(8.8)
 $$
where $  L_j =L_j(z), \ M_j = M_j(z), \ z \in (0,\infty), \ j=1,2,\ldots, l $ are
{\it slowly varying} simultaneously as $ z \to 0  $  and as $ z \to \infty $ continuous
positive functions, $ \vec{p} = \{ p_j \}, \ \vec{q} = \{ q_j \}, $
$$
 p_j \ge p_{j,0} \stackrel{def}{=} \frac{m_j}{m_j-\beta_j}, \ 1 \le q_j <
 q_{j,0} \stackrel{def}{=} \frac{m_j}{\alpha_j},
$$
$ 0 \le \alpha_j, \beta_j < 1. $ \par

\vspace{3mm}

{\bf Theorem 8.2.}\par

\vspace{3mm}

{\bf A.} It the inequality (8.8) there holds for some non-zero function
$ f, \ f \in S(R^n), $ then

$$
\frac{\beta_j - \alpha_j}{m_j} = 1 - \frac{1}{p_j} -  \frac{1}{q_j}, \eqno(8.9)
$$

$$
L_j(1/\lambda) \asymp M_j(\lambda). \eqno(8.10)
$$

\vspace{2mm}

{\bf B.} If the equations (8.9) and (8.10) are satisfied, then

 $$
 C_5(L,M; \vec{\alpha},\vec{\beta,\vec{m}})  \prod_{j=1}^l
 \left[ \frac{p_j}{(p_j-m_j/(m_j-\beta_j))}\right]^{(\alpha_j+\beta_j)/m_j}  \le
K_{L,M;n,\vec{\alpha},\vec{\beta}}(\vec{p}) \cdot
 $$

$$
C_6(L,M; \vec{\alpha},\vec{\beta},\vec{m}) \prod_{j=1}^l
\left[\frac{p_j}{(p_j-m_j/(m_j-\beta_j))} \right]^{\max(1,(\alpha_j+\beta_j)/m_j)}.
 \eqno(8.11)
$$
(anisotropic version).

\vspace{3mm}

 Evidently, if

$$
\frac{\beta_j - \alpha_j}{m_j} = \const,  \ j = 1,2,\ldots,l;
$$
then  for the values

$$
p > \max_j \frac{m_j}{m_j-\beta_j}
$$
there holds
 $$
 | \ |y|^{- \vec{\alpha}} F[f](y) \ \prod_{j=1}^l M_j(|y_j|)|_q
 \le K_{L,M;n,\vec{\alpha},\vec{\beta}}(p) \ \cdot
 | \ |x|^{\vec{\beta}} \ f(x) \ \prod_{j=1}^l L_j(|x_j|)|_p, \eqno(8.12)
 $$

$$
\frac{\beta_j - \alpha_j}{m_j} = 1 - \frac{1}{p} - \frac{1}{q},
$$

 $$
 C_5(L,M; \vec{\alpha},\vec{\beta,\vec{m}})  \prod_{j=1}^l
 \left[ \frac{p}{(p-m_j/(m_j-\beta_j))}\right]^{(\alpha_j+\beta_j)/m_j}  \le
K_{L,M;n,\vec{\alpha},\vec{\beta}}(p) \cdot
 $$

$$
C_6(L,M; \vec{\alpha},\vec{\beta},\vec{m}) \prod_{j=1}^l
\left[\frac{p}{(p-m_j/(m_j-\beta_j))} \right]^{\max(1,(\alpha_j+\beta_j)/m_j)}.
\eqno(8.13)
$$
(isotropic version).\par

\vspace{3mm}
\section{Concluding remarks.}

\vspace{3mm}

{\bf 1.} We introduce the so-called {\it weight} $ L_p $ norm as follows. Let
$  \mu $ be arbitrary constant. Denote for the measurable function
$ g: R^n \to R $

$$
|g|_{p,\mu} = \left[\int_{R^n} |g(x)|^p \ |x|^{\mu} \ dx  \right]^{1/p}.
$$

 Let the parameters $ \alpha,\beta,\lambda,\mu; p,q $  be such that
 $ n-\lambda > \alpha, \ q_0:=(n-\lambda)/\alpha >1, $
 $$
 \alpha,\beta \in (0,1), \ \lambda,\mu \ge 0,
 \ p > p_0 := (n+\mu)/(n-\beta), q \in [1,q_0),
 $$

$$
\beta-\alpha = n - n \left(\frac{1}{p+\mu/\beta}  + \frac{1}{q+\lambda/\alpha} \right).
$$
 The proposition of theorem  4.1 may be rewritten after change of parameters as follows.
Let us denote

$$
K_{\lambda,\mu}(p) = \sup_{f \ne 0, \ \in L_{p,\mu}}
\left[ \frac{ | \ |x|^{-\alpha} \ F[f]|_{q,-\lambda}}
{| \ |x|^{\beta} \ |f| \ |_{p,\mu}}\right];
$$
then

$$
C_1(\alpha,\beta,\lambda,\mu)[p/(p-p_0)]^{1/q_0+1/p_0 } \le K_{\lambda,\mu}(p) \le
$$

$$
C_2(\alpha,\beta,\lambda,\mu)[p/(p-p_0)]^{\max(1, 1/q_0 + 1/p_0)}, \ p > p_0.
$$

\vspace{3mm}

{\bf 2.} We introduce and calculate in this pilcrow the multidimensional Boyd's indices
for anisotropic Grand Lebesgue Spaces, which play a very important role in the theory
of Fourier series \cite{Bennet}, chapter 6.7.\par
 In detail, let $ X = R^1_+ \times R^1_+ $ with
ordinary Lebesgue measure. Let $  G_D \psi, \ \psi: (a,b)\times (c,d) \to R_+$
be some function from the set $  \Psi_D, \ D = (a,b)\times (c,d), 1 \le a < b < \infty,
1 \le c < d \le \infty. $ We introduce  the multidimensional Boyd's indices as follows.
Denote for the values $ s,t > 0 $ the multidimensional (in our case- two dimensional)
dilation operator
$$
\Delta_{s,t}[f](x,y) = f(x/s,y/t),
$$
and define

$$
\overline{\alpha}(G\psi_D) \stackrel{def}{=}
\lim_{s \to \infty} \frac{\log \left[||\Delta_{s,t}||G\psi_D \right]}{|\log s|},
$$

$$
\underline{\alpha}(G\psi_D) \stackrel{def}{=}
\lim_{s \to 0+} \frac{\log \left[ ||\Delta_{s,t}||G\psi_D \right]}{|\log s|},
$$

$$
\overline{\beta}(G\psi_D) \stackrel{def}{=}
\lim_{t \to \infty} \frac{\log  \left[||\Delta_{s,t}||G\psi_D \right]}{|\log t|},
$$

$$
\underline{\beta}(G\psi_D) \stackrel{def}{=}
\lim_{t \to 0+} \frac{ \left[\log ||\Delta_{s,t}||G\psi_D \right]}{|\log t|}.
$$
 We conclude analogously to the article \cite{Ostrovsky3}:

 $$
\overline{\alpha}(G\psi_D) = 1/a, \ \underline{\alpha}(G\psi_D)  = 1/b,
 $$

 $$
\overline{\beta}(G\psi_D) = 1/c, \ \underline{\beta}(G\psi_D)  = 1/d.
 $$

\vspace{3mm}

{\bf 3.} The "dual" version of PBO inequality may be formulated as follows:

$$
| \ |x|^{-\alpha} \ f(x) \ |_q \le K_{PBO}(p) \ | \ |y|^{\beta} \ F[f](y) \ |_p,
$$
where as before

$$
1 \le q < n/\alpha, \ p > n/(n-\beta), \ \beta - \alpha = n(1-1/p - 1/q).
$$

\vspace{3mm}

{\bf 4.} In the terms of anisotropic norms may be estimated the $ L_p \to L_q $ 
norm of the integral operators of a view

$$
W[f](x) = \int_Y K(x,y) \ f(y) \ \nu(dy), \ x \in X.
$$
 Namely, it is proved in \cite{Jorgens1}, p. 272 that 
 
$$
|W|(L_p \to L_q) \le |K|_{q_1,p}, \ q_1 = q/(q-1), \ q > 1, p \in (1,\infty).
$$

\end{document}